\documentclass[a4paper]{ifacconf}
\DeclareUnicodeCharacter{02B9}{'}
\DeclareUnicodeCharacter{0308}{"}

\usepackage{graphicx}      
\usepackage{natbib}        
\usepackage{mathptmx}
\usepackage{amssymb}
\usepackage{amsmath}
\usepackage{caption}
\usepackage{mathdots}
\usepackage{bm}

\usepackage{algorithmic}
\usepackage{algorithm}
\usepackage[ruled,vlined,algo2e]{algorithm2e}

\usepackage[all,cmtip]{xy}
\usepackage{tikz}
\usetikzlibrary{automata,arrows, positioning, shapes.geometric, shadows.blur, fit} 
\begin{document}
\begin{frontmatter}

\title{Stochastically Structured Reservoir Computers for Financial and Economic System Identification}


\author[First]{Lendy Banegas}\hspace{2pc} 
\author[Second]{Fredy Vides}

\address[First]{Department of Statistics and Research, CNBS, Honduras (e-mails: lendy.banegas@cnbs.gob.hn).}


\address[Second]{Center for Innovation in Scientific Computing, School of Mathematics and Computer Science, Faculty of Science, UNAH, Honduras (e-mail:  fredy.vides@unah.edu.hn).}

\begin{abstract}                
This paper introduces a methodology for identifying and simulating financial and economic systems using stochastically structured reservoir computers (SSRCs). The framework combines structure-preserving embeddings with graph-informed coupling matrices to model inter-agent dynamics while enhancing interpretability. A constrained optimization scheme guarantees compliance with both stochastic and structural constraints. Two empirical case studies, a nonlinear stochastic dynamic model and regional inflation network dynamics, demonstrate the effectiveness of the approach in capturing complex nonlinear patterns and enabling interpretable predictive analysis under uncertainty.
\end{abstract}


\end{frontmatter}

\section{Introduction}

This paper introduces a methodology for identifying and simulating financial and economic systems using stochastically structured reservoir computers (SSRCs). The approach builds on reservoir computing while constraining the output coupling matrices through relational graphs that reflect observed inter-agent interactions. Combined with stochastic embeddings, this design preserves structural properties and improves interpretability in system identification.

The methodology is illustrated through two empirical case studies: (i) a nonlinear stochastic dynamic model that highlights how SSRCs capture probabilistic dynamics and feedback in adaptive systems, and (ii) regional inflation dynamics across the CAPARD region (Central America, Panama, and the Dominican Republic), the United States, and China. The latter highlights the role of U.S. monetary policy, via changes in the federal funds rate, in shaping domestic and regional inflation during and after the COVID-19 pandemic.

Relative to recent advances in reservoir computing (Yan et al., 2024; Ehlers et al., 2025), this work makes three contributions. First, it integrates stochastic matrix structures into the framework, enhancing its capacity to capture nonlinear and uncertain economic behaviors. Second, it develops a constrained optimization scheme that balances structural compliance with computational efficiency. Third, by focusing on financial and economic applications, it extends the scope of reservoir computing to domains where interpretable, data-driven modeling is essential for policy and decision making.

The remainder of the paper is organized as follows. Section \S2 introduces the preliminaries and notation. Section \S3 presents the proposed methodology, followed by algorithms and interpretive discussion in Section \S4. Section \S5 reports simulation results for the two case studies, and Section \S6 concludes with future directions.

\section{Preliminaries and Notation}

In this study we will identify the space of vectors $\mathbb{R}^n$ with (column) matrices in $\mathbb{R}^{n\times 1}$. The symbol $\mathbb{R}_{0}^{+}$ will be used to denote the positive real numbers including zero. We will write $\mathbf{1}_n$ to denote the (column) vector in $\mathbb{R}^n$ with all its entries equal to one.

For any vector $x\in \mathbb{R}^n$, we will write $x[j]$ to denote the jth component of $x$.

Given $\delta>0$, we will denote by $H_\delta$ the function defined by the expression 
\begin{align*}
H_\delta(x)=\left\{
\begin{array}{ll}
1, & x>\delta\\
0,& x\leq \delta
\end{array}
\right.,
\end{align*}
for any $x\in \mathbb{R}$.

Let \( X \in \mathbb{R}^{n \times n} \) with \( x_{ij} \) as its \((i,j)\)-th element. The vectorization of \( X \), denoted \( \operatorname{vec}(X) \), stacks its columns into:
\[
\operatorname{vec}(X) := 
\begin{bmatrix}
x_{11} & x_{21} & \dots & x_{n1} & x_{12} & \dots & x_{nn}
\end{bmatrix}^T \in \mathbb{R}^{n^2}.
\]

The set \(\mathbb{S}_{m,n}(\mathbb{R})\) defined as:
\[
\mathbb{S}_{m,n}(\mathbb{R}) := \left\{ A \in (\mathbb{R}^+_0)^{m \times n} \ \middle| \ \mathbf{1}_m^\top A = \mathbf{1}_n^\top  \right\}
\]
represents the class of stochastic (column-stochastic) matrices in $\mathbb{R}^{m\times n}$. We will say that a vector $x\in(\mathbb{R}_0^+)^n$ is substochastic if $\mathbf{1}_n^\top x\leq 1$. The symbol $e_{j,k}(m,n)$ will denote the matrix in $\mathbb{R}^{m\times n}$ with $j,k$ entry equal to one and with zeros elsewhere.

For any integer $p>0$ and any matrix $X\in \mathbb{R}^{m\times n}$, we will write $X^{\otimes p}$ to denote the operation determined by the following expression.
\[
X^{\otimes p} =\left\{
\begin{array}{ll}
X&, p=1\\
X\otimes X^{\otimes (p-1)}&, p\geq 2
\end{array}
\right.
\]

\section{Methodology}

\label{sec:methods}

The dynamic models considered in this study are described by discrete-time systems of the form
\begin{align}
y(t) &= W\,\eth_p(x(t)) + e_t,
\label{eq:generic_control_system_1}
\end{align}
where, for each $t > 0$, the vectors $x(t)$ and $y(t)$ are stochastic vectors of appropriate dimensions, and the term $e_t := y(t) - W\,\eth_p(x(t))$ denotes the residual component.

The models in \eqref{eq:generic_control_system_1} correspond to what we refer to in this study as stochastically structured reservoir computers. Their architecture generalizes the framework introduced in \cite{modelos_reservorios}, with $W$ constrained to be a column-stochastic matrix and with the mapping $\eth_p$ satisfying the property that $\eth_p(z)$ is stochastic whenever $z$ is a stochastic vector of appropriate dimension.

For this study we focus on a specific family of mappings \( \eth_p \).
Given a vector \( s = [s_1~\cdots~s_{p+1}]^\top \in \{0,1\}^{p+1} \), let
\( \tilde{\eth}_{s,p}(x) : \mathbb{R}^n \to \mathbb{R}^{d_p(n)} \) denote the map
\[
\tilde{\eth}_{s,p}(x) := \frac{1}{m}
\begin{bmatrix}
s_1\, x \\
s_2\, x^{\otimes 2} \\
\vdots \\
s_p\, x^{\otimes p} \\
s_{p+1}
\end{bmatrix},
\]
where \( d_p(n) = \frac{n(n^{p}-1)}{n-1} + 1 \), and \( m = \sum_{i=1}^{p+1} s_i \)
denotes the number of active blocks. The parameter $m$, defined earlier, denotes the number of block terms considered. The selection of terms may be guided by exploratory analysis of prior signals or by domain knowledge. The integer $p$ will be referred to as the \emph{maximum order} of the embedding map $\tilde{\eth}_{s,p}$. From here on, and for computational purposes, $\tilde{\eth}_{s,p}(x)$ will denote the expression obtained after eliminating the zero blocks from its formal definition. Under this consideration, the resulting map $\tilde{\eth}_{s,p}$ can now be seen as a mapping  $\tilde{\eth}_{s,p}:\mathbb{R}^n\to \mathbb{R}^{d_{s,p}(n)}$, for $d_{s,p}(n):=s_{p+1}+\sum_{k=1}^{p}s_k n^k$.\\

\begin{lem}\label{lem:structured_embedding}
The embedding map $\tilde{\eth}_{s,p}$ preserves stochastic vectors.
\end{lem}
\begin{pf}
It can be seen that for any $x,y\in \mathbb{S}_{n,1}(\mathbb{R})$, if we define $z=x\otimes y$, then $z\in \mathbb{R}^{n^2}$ and we will have that:
\begin{align*}
    \mathbf{1}_{n^2}^\top z&=\sum_{j=1}^n x[j]\sum_{k=1}^n y[k]=\sum_{j=1}^n x[j]\cdot 1 = \sum_{j=1}^n x[j] = 1.
\end{align*}
This implies that for any $p$, $x^{\otimes p}$ is stochastic, and by the definition of $\tilde{\eth}_{s,p}$, we will have that the sum of the entries of $\eth_{s,p}(x)$ equals $m/m=1$. Consequently, the sum of the entries of $\tilde{\eth}_{s,p}(x)$ equals 1. This completes the proof.\\
\end{pf}
From here on, a map $\tilde{\eth}_{s,p}$ will be called a stochastic $p$-embedding.\\

\begin{lem}\label{lem:Compression_existence}
Given positive integers $n,p$ and $s\in \{0,1\}^{p+1}$. There are an integer $0<\rho_p(n)<d_{s,p}(n)$ and a sparse matrix $R_{s,p}(n)\in \mathbb{R}^{\rho_p(n)\times d_{s,p}(n)}$ with $d_{s,p}(n)$ nonzero entries, such that $R_{s,p}(n)\: \tilde{\eth}_{s,p}(x)$ is stochastic and has the least number of non-redundant words (monomial terms) for any $x\in \mathbb{R}^{n}$.
\end{lem}
\begin{pf}
Let $n,p$ be positive integers. Given $s\in \{0,1\}^{p+1}$, consider the structured embedding map $\tilde{\eth}_{s,p}:\mathbb{R}^{n} \to \mathbb{R}^{d_{s,p}(n)}$, where $d_{s,p}(n):=s_{p+1}+\sum_{k=1}^{p}s_k n^k$ corresponds to the total number of distinct tensor monomials up to degree $p$, plus a constant term.

We aim to construct a sparse matrix $R_{s,p}(n) \in \mathbb{R}^{\rho_p(n) \times d_{s,p}(n)}$ that maps the embedding $\tilde{\eth}_{s,p}(x)$ to a stochastic vector of reduced dimension, while preserving all non-redundant monomial terms.

Let us start by defining the matrix $R\in \mathbb{R}^{1\times d_{s,p}(n)}$, with $1$ in its $R_{11}$ entry and with all other entries equal to zero.

For each $2\leq j\leq d_{s,p}(n)$, let us consider the indices $j=k_1(j)< k_2(j)<\cdots <k_{n_j}(j)<d_{s,p}(n)$ that correspond to the same monomial in $\tilde{\eth}_{s,p}(x)$, let us define the matrix $R_0\in \mathbb{R}^{1\times d_{s,p}(n)}$ with $1$ in its $R_{1k_l(j)}$ entries for $1\leq l\leq n_j$, and with all other entries equal to zero. Let us now define the augmented matrix 
$$
R:=\begin{bmatrix}
    R\\
    R_0
\end{bmatrix}
$$
Finally, update the matrix $R$, using the operation:
$$
R:=\begin{bmatrix}
    R\\
    R'
\end{bmatrix}
$$
where $R'\in \mathbb{R}^{1\times d_{s,p}(n)}$ is the matrix with entry $R'_{1d_{s,p}(n)}$ equal to $1$, and with all other entries equal to $0$. 

Let us set $R_{s,p}(n):=R$. It can be seen by the way $R$ has been constructed, that the operation $R_{s,p}(n)\tilde{\eth}_{s,p}(x)$ assigns each group of replicates in $\tilde{\eth}_p(x)$ to a single representative coordinate of $R_p(n)\tilde{\eth}_p(x)$. This selection is performed by adding over the redundant entries and replacing this sum into one suitable coordinate. Because of this, and as a consequence of Lemma \ref{lem:structured_embedding}, it is clear that $R_{s,p}(n) \tilde{\eth}_{s,p}(x)$ is stochastic. Let us set $\rho_{p}(n)$ as the number of rows of $R$. This completes the proof.
\end{pf}

Using the matrix \( R_{s,p}(n) \) described in Lemma~\ref{lem:Compression_existence}, we define the reduced stochastic embedding \( \tilde{\eth}_{s,p,r} \) as  
\[
\tilde{\eth}_{s,p,r}(x) := R_{s,p}(n)\tilde{\eth}_{s,p}(x),
\]  
for any \( x \in \mathbb{S}_{n,1}(\mathbb{R}) \). This transformation yields a compressed representation of the full embedding, preserving only the non-redundant monomial components while ensuring stochasticity.

For formal implementation purposes, the operation determined by the expression \( R_{s,p}(n)\tilde{\eth}_{s,p}(x) \) can be equivalently represented by a vector composed of the distinct monomials of degree less than or equal to \( p \), each multiplied by the appropriate integer scaling factor that accounts for its multiplicity in the original tensor expansion.

In this context, a \emph{structured reservoir computer} refers to a model of the form \eqref{eq:generic_control_system_1}, with $\eth_p$ denoting either a reduced stochastic embedding $\tilde{\eth}_{s,p,r}$ or a function of its terms, as illustrated in Section~\S5. The output coupling matrix $W$ in \eqref{eq:generic_control_system_1} is stochastic and its structure is determined by a relational graph $\mathcal{G}_{\mathcal{S}}=(V_\mathcal{S},E_{\mathcal{S}})$, representing observed or expected interactions among dynamic agents from $t$ to $t+1$ within the training time frame $\{0,1,\ldots,T\}$, for some $T>0$.

Let us consider any system that describes the dynamic interation of $n>0$ agents, then $V_{\mathcal{S}}=\{1,\ldots,n\}$, and if $\mathbb{B}_{\mathcal{S}}(n)$ denotes the set
\begin{equation}
\mathbb{B}_{\mathcal{S}}(m,n):= \left\{e_{j,k}(m,n)\left|(j,k)\in E_{\mathcal{S}}\right.\right\},
    \label{eq:structure_matrix_set}
\end{equation}
 then the existence of a $\mathbb{B}_{\mathcal{S}}(m,n)$-structured stochastic output coupling matrix $W$ that satisfies \eqref{eq:generic_control_system_1}, is a consequence of the following theorem.\\

\begin{thm}
    The coupling matrices $W$ in equation \eqref{eq:generic_control_system_1} can be approximately identified in the matrix set $\left(\mathrm{span}\:\mathbb{B}_{\mathcal{S}}(m,n)\right)\cap \mathbb{S}_{m,n}(\mathbb{R})$.
\end{thm}

\begin{pf}
Let $\{x(t) \in \mathbb{R}^n \mid t = 0, 1, \ldots, T\}$ and $\{y(t) \in \mathbb{R}^m \mid t = 0, 1, \ldots, T\}$ denote a stochastic vector time series corresponding to the system's evolution. We begin by applying a function $\eth_p$ of a reduced embedding transformation $\tilde{\eth}_{s,p,r} := R_{s,p}(n) \tilde{\eth}_{s,p}$ to the observed states, where $p$ is a prescribed tensor order and $R_{s,p}(n)$ is the structure-preserving compression map guaranteed by Lemma~\ref{lem:Compression_existence}. This yields a compressed representation of each state, which we collect into the input data matrix
\[
X_0(T) := \begin{bmatrix}
| &  & |\\
\eth_{p}(x(0)) & \cdots & \eth_{p}(x(T-1)) \\
| &  & |
\end{bmatrix},
\]
while the elements $y(t)$ are stored in
\[
X_1(T) := \begin{bmatrix}
| &  & |\\
y(0) & \cdots & y(T-1) \\
| &  & |
\end{bmatrix}.
\]

To ensure that the identified coupling matrix respects the structural constraints of the system, we consider a structured dictionary $\mathbb{B}_{\mathcal{S}}(m,n) = \{S_1, \ldots, S_q\} \subset \mathbb{R}^{m \times n}$, derived from the graph $\mathcal{G}_\mathcal{S}$ that encodes allowable interactions. Each $S_j$ acts as a basis component capturing a localized or interpretable mode of coupling. We construct a matrix $X_0$ whose columns correspond to the vectorized actions of these basis elements on the embedded data:
\[
X_0 := \begin{bmatrix}
\text{vec}(S_1 X_0(T)) & \cdots & \text{vec}(S_q X_0(T))
\end{bmatrix},
\]
and define the target vector as $X_1 := \text{vec}(X_1(T))$.

The identification task is then posed as the following structured nonnegative least squares problem:
\[
\mathbf{a}:=\arg\min_{\hat{\mathbf{a}}\in (\mathbb{R}_0^+)^q}\left\|\begin{bmatrix}
X_0^\top X_0 \\
C
\end{bmatrix} \hat{\mathbf{a}}
-
\begin{bmatrix}
X_0^\top X_1 \\
\mathbf{1}_p
\end{bmatrix}\right\|_2,
\]
where $\|\cdot\|_2$ denotes the Euclidean vector norm, and $C$ is a constraint matrix chosen to enforce the membership of the resulting linear combination $\hat{W} := \sum_{j=1}^q \mathbf{a}[j] S_j$ in the structured subset $\mathbb{S}_{m,n}(\mathbb{R})$. Such constraints may include nonnegativity, block sparsity, or stochasticity properties, depending on the system's prior assumptions.

As discussed in \cite{BOUTSIDIS2009760}, the above system corresponds to a convex quadratic optimization problem with linear constraints (e.g., a Nonnegative Least Squares problem), and is therefore solvable in polynomial time up to arbitrary precision.

Since $\hat{W}$ is constructed as a linear combination of elements in $\mathbb{B}_{\mathcal{S}}(m,n)$, it lies in $\mathrm{span}\:\mathbb{B}_{\mathcal{S}}(m,n)$ by definition. The enforcement of structural constraints via $C$ ensures that $\hat{W} \in \mathbb{S}_{m,n}(\mathbb{R})$. Hence,
\[
\hat{W} \in \left( \mathrm{span}\:\mathbb{B}_{\mathcal{S}}(m,n) \right) \cap \mathbb{S}_{m,n}(\mathbb{R}),
\]
as claimed. This completes the proof.
\end{pf}

\section{Algorithms}
\label{sec:algorithm}

In this section, we focus on the applications of the structured matrix approximation methods presented in \S\ref{sec:methods}, to reservoir computer models identification for stochastically structured dynamical systems. More specifically, we propose a prototypical algorithm for general purpose stochastically structured system identification, that is described by Algorithm \ref{alg:main_AutoRegressor_alg_2}, and that is based on the structured least squares solver described by Algorithm \ref{alg:main_SLMESolver_alg_1}.

\begin{algorithm2e}
\caption{{\bf SSRC Model}: SSRC model identification}
\label{alg:main_AutoRegressor_alg_2}
\SetAlgoLined
 \KwData{$\Sigma_{T}(x)=\{x(t)\}_{t=1}^{T}\subset \mathbb{R}^n$, $\Sigma_{T}(y)=\{y(t)\}_{t=1}^{T}\subset \mathbb{R}^m, \mathbb{B}_{\mathcal {S}}(m,n)\subset \mathbb{R}^{m\times q}$}
  \KwResult{$\hat{W},\eth_{p}$}
\begin{itemize}
\item [1:] Choose or identify $s\in \{0,1\}^{p+1}$.
\item[2:] Compute $\eth_p$ as either the reduced embedding $\tilde{\eth}_{s,p,r}$ or a function of it, using the map $R_{s,p}(n)$ from Lemma~\ref{lem:Compression_existence}.

\item[3:] Compute matrices:
\begin{align*}
X_0(T)&:=R_p(n)\begin{bmatrix}
| &  & |\\
    \eth_p\left(x(0)\right) & \ldots & \eth_p\left(x(T-1)\right)\\
    | &  & |
\end{bmatrix} \label{eq:training_data_matrix}\\
X_1(T)&:=\begin{bmatrix}
| &  & |\\
    y(0) & \ldots & y(T-1)\\
    | &  & |
\end{bmatrix}
\end{align*}
\item[4:] Set:
\begin{align*}
    X_0& = \begin{bmatrix}
    | &  & | \\
        vec(S_1 X_0(T)) & \cdots & vec(S_q X_0(T))\\
        | &  & |
    \end{bmatrix}\\
    X_1& = vec(X_1(T))
\end{align*}
for $\mathbb{B}_{\mathcal{S}}(m,n):=\{S_1,\ldots,S_q\}$ 
\item[5:] Approximately solve the structured nonnegative least squares problem:
\begin{equation*}
\mathbf{a}:=\arg\min_{\hat{\mathbf{a}}\in (\mathbb{R}_0^+)^q}\left\|\begin{bmatrix}
X_0^\top X_0 \\
C
\end{bmatrix} \hat{\mathbf{a}}
-
\begin{bmatrix}
X_0^\top X_1 \\
\mathbf{1}_p
\end{bmatrix}\right\|_2
\end{equation*}
using Algorithm \ref{alg:main_SLMESolver_alg_1}

\item[6:] Set: 
\[
\hat{W}:=\sum_{j=1}^q \mathbf{a}[j]S_j
\]
\end{itemize}
\KwRet{$\hat{W},\eth_{p}$}
\end{algorithm2e}

\begin{algorithm2e}
\caption{{\bf SSSLSSolver}: Sparse stochastically structured linear least squares solver algorithm}
\label{alg:main_SLMESolver_alg_1}
\SetAlgoLined
\KwData{$A\in \mathbb{R}^{m\times n}$, $y\in \mathbb{R}^{m}$, $\delta>0$, $N\in \mathbb{Z}^+$, $\varepsilon>0$}
\KwResult{$x$}
\BlankLine
Set $x_0=\arg\min_{\tilde{c}\in (\mathbb{R}^+_0)^{n}}\|A\tilde{c}-y\|_2$\;
    Set $K=1$\;
    Set $\mathrm{error}=1+\delta$\;
    Set $c=x_0$\;
    Set $\hat{c}=\begin{bmatrix}
    \hat{c}_1 & \cdots & \hat{c}_n
    \end{bmatrix}^\top=\begin{bmatrix}
    |\hat{e}_{1,n}^\ast c| & \cdots & |\hat{e}_{n,n}^\ast c|
    \end{bmatrix}^\top$\;
    Compute permutation $\sigma:\{1,\ldots,n\}\to \{1,\ldots,n\}$ such that $\hat{c}_{\sigma(1)}\geq \hat{c}_{\sigma(2)}\geq \cdots \geq \hat{c}_{\sigma(n)}$\;
    Set $N_0=\max\left\{\sum_{j=1}^n H_\varepsilon\left(\hat{c}_{\sigma(j)}\right),1\right\}$\;
    \While{$K\leq N$ $\And$ $\mathrm{error}>\delta$}{
        Set $x=\mathbf{0}_{n}$\;
        Set $A_0=\sum_{j=1}^{N_0} A\hat{e}_{\sigma(j),n}\hat{e}_{j,N_0}^\ast$\;
        Solve $c=\arg\min_{\tilde{c}\in (\mathbb{R}^+_0)^{N_0}}\|A_0\tilde{c}-\hat{Y}\hat{e}_{j,p}\|_2$\;
        \For{$k=1,\ldots,N_0$}{
            Set $x_{\sigma(k)}=\hat{e}_{k,N_0}^\ast c$\;
        }
        Set $\mathrm{error}=\|x-x_0\|_2$\;
        Set $x_0=x$\;
        Set $\hat{c}=\begin{bmatrix}
        \hat{c}_1 & \cdots & \hat{c}_n
        \end{bmatrix}^\top=\begin{bmatrix}
        |\hat{e}_{1,n}^\ast x| & \cdots & |\hat{e}_{n,n}^\ast x|
        \end{bmatrix}^\top$\;
        Compute permutation $\sigma:\{1,\ldots,n\}\to \{1,\ldots,n\}$ such that $\hat{c}_{\sigma(1)}\geq \hat{c}_{\sigma(2)}\geq \cdots \geq \hat{c}_{\sigma(n)}$\;
        Set $N_0=\max\left\{\sum_{j=1}^n H_\varepsilon\left(\hat{c}_{\sigma(j)}\right),1\right\}$\;
        Set $K=K+1$\;
    }
\KwRet{$x$}
\end{algorithm2e}

\subsection{Interpretation}

From a control and optimization standpoint, the proposed framework formalizes systemic adaptation through structure-constrained embeddings and graph-restricted coupling matrices. The stochastic embedding \( \tilde{\eth}_{p,r}(x) \) acts as a dimensional and structural regularizer, ensuring that the reconstructed representation remains consistent with the probabilistic and relational constraints of the system.

Accordingly, the dynamics in~\eqref{eq:generic_control_system_1} perform more than state prediction: the embedding encodes feasible transformations within the stochastic manifold, while the optimization stage identifies the coupling matrix \( W \) inside the admissible structural subspace \( (\mathrm{span}\,\mathbb{B}_\mathcal{S}(m,n)) \cap \mathbb{S}_{m,n}(\mathbb{R}) \). This joint mechanism promotes interpretability and stability in the identified model.

Thus, the SSRC algorithm operates as a structure-preserving solver that unifies identification, constraint enforcement, and adaptive regularization, providing a transparent link between algorithmic operations and system-level behavior.

\medskip

\noindent\textbf{Remark.}
Models of the form~\eqref{eq:generic_control_system_1} also allow the identification of autoregressive systems with control inputs. This is achieved by replacing \( y(t) \) with \( x(t+1) \) and redefining the embedding input as
\[
x(t)\;\mapsto\;
\begin{bmatrix}
x(t)\\[2pt]
u(t)
\end{bmatrix},
\]
where \( u(t) \) denotes an exogenous or control signal.

\section{Computational simulations}
\label{sec:results}

The dynamic models considered in this section are determined by particular representations of the dynamic models \eqref{eq:generic_control_system_1}.

\subsection{Modeling Preference Evolution in Investment Coalitions with SSRCs}

This simulation introduces a stochastic dynamic model of resource allocation among agents who interact strategically and adaptively to maximize their participation benefits in the selection of an investment account within a given portfolio. Rather than advancing normative claims about the desirability of specific participation levels, the analysis focuses on the evolving trajectories of allocation and the emergent patterns over time. The model considered is given by  
\begin{equation}
p(t+1) = W\,\eth_2(p(t)),
\label{eq:example_model_1}
\end{equation}  
where $\eth_2(p(t))$ denotes a structured embedding of reduced second-degree monomials extracted from $\tilde{\eth}_{(1,1,0),2}(p(t))$, defined as  
\[
\eth_2(p(t)) =
\begin{bmatrix}
    0.9\,p(t)\\
    0.1\,p(t)[1]\,p(t)[2]\\
    0.1\,p(t)[3]\,p(t)[5]\\
    0.1\,(1 - p(t)[1]\,p(t)[2] - p(t)[3]\,p(t)[5])
\end{bmatrix},
\]  
with $p(t)[j]$ denoting the probability of choosing account $j$ in the portfolio at time $t$.

The structural configuration of the maps in \eqref{eq:example_model_1} allows us to interpret the system as a nonlinear dynamic model of the form:  
\begin{equation}
p(t+1) = 0.9\,W_1\,p(t) + 0.1\,W_2\,u(p(t)),
\label{eq:closed_loop_rep_model_1}
\end{equation}  
where $u(p(t))$ is a stochastic vector of reduced second-degree monomials, and $W_1, W_2$ are stochastic matrices that represent state transitions and feedback interactions, respectively:  
\begin{align*}
    W_1 &:= \begin{bmatrix}
        1 & 0.45 & 0 & 0 & 0\\
        0 & 0 & 0.45 & 0 & 0\\
        0 & 0.55 & 0 & 0.45 & 0\\
        0 & 0 & 0.55 & 0 & 0\\
        0 & 0 & 0 & 0.55 & 1
    \end{bmatrix}, \quad
    W_2 := \begin{bmatrix}
        0.8 & 0.1 & 0\\
        0 & 0 & 0\\
        0 & 0 & 0\\
        0 & 0 & 0\\
        0.2 & 0.9 & 1
    \end{bmatrix}.
\end{align*}

\begin{figure}[H]
    \centering
    \includegraphics[width=0.95\linewidth]{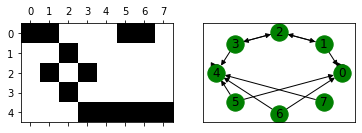}
    \caption{Structured coupling matrix (left) and relational graph (right) of the investment preference evolution model.}
    \label{fig:coupling-matrix}
\end{figure}

The structured output coupling matrix $\hat{W}$ identified from the model \eqref{eq:example_model_1}, using synthetic signals generated by iterating on it, is given by  
\[
\hat{W} := \begin{bmatrix}
    1 & 0.45 & 0 & 0 & 0 & 0.8 & 0.1 & 0\\
    0 & 0 & 0.45 & 0 & 0 & 0 & 0 & 0\\
    0 & 0.55 & 0 & 0.45 & 0 & 0 & 0 & 0\\
    0 & 0 & 0.55 & 0 & 0 & 0 & 0 & 0\\
    0 & 0 & 0 & 0.55 & 1 & 0.2 & 0.9 & 1
\end{bmatrix},
\]
with its sparsity pattern shown in Figure~\ref{fig:coupling-matrix}, together with the associated relational graph.

\begin{figure}[H]
    \centering
    \includegraphics[width=0.9\linewidth]{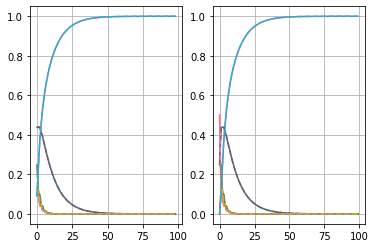}
    \caption{Estimated short-term evolution (left) and long-term evolution (right) of investment preference allocations.}
    \label{fig:signals-id}
\end{figure}

The results of the identification process based on models~\eqref{eq:example_model_1} and~\eqref{eq:closed_loop_rep_model_1} are shown in Figure~\ref{fig:signals-id}, illustrating the short- and long-term evolution of investment preference allocations.

This example shows that SSRCs can model stochastically embedded substochastic systems. By Lemma~\ref{lem:structured_embedding}, any sub-vector extracted from tensor product blocks under a stochastic embedding defines a substochastic vector. This property will also be used in the next example to identify relational graphs among economic variables.

Thanks to their flexible structure, models of the form~\eqref{eq:closed_loop_rep_model_1} apply across social, economic, financial, ecological, and symbolic domains. Although they do not impose a predefined causal structure, the results raise questions about the conditions under which resource concentration may be linked to broader dynamics such as stratification, exclusion, or power asymmetries.

\subsection{Regional Inflation Network Dynamics} 

Several studies have highlighted the significant influence of U.S. monetary policy on the monetary frameworks and inflation rates of other countries \citep{Tenkovskaya2023, Azad2022, Carella2024}. Furthermore, the existence of cross-country inflation transmission dynamics has also been documented \citep{Budova2023, Iraheta2008, Liu2015}.

Following the COVID-19 recession, many countries, including the United States, adopted expansionary monetary policies to mitigate the economic downturn, with the federal funds rate being a key tool \citep{Feldkircher2021}. According to \cite{Swanson2024}, changes in the federal funds rate significantly impact production and prices, highlighting short-term interest rates as the most powerful tool central banks have to influence the economy.

\begin{figure}[!ht]
    \centering
    \includegraphics[scale=.55]{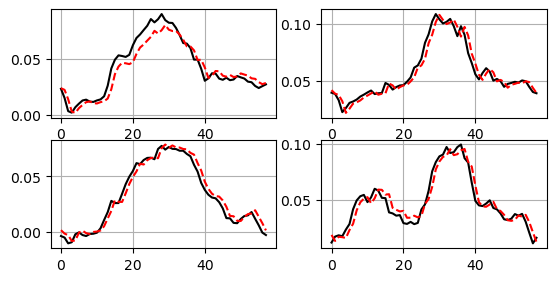}
    \caption{US inflation signal (top left). Honduran inflation signal (top right). Salvadoran inflation signal (bottom left). Guatemalan inflation signal (bottom right).}
    \label{fig:Eco_Signals}
\end{figure}

Beginning in 2022, the Federal Reserve adjusted its strategy in response to external factors and rising inflation, leading to higher federal rates and tighter financial conditions \citep{Alekseievska2024}. Most of the CAPARD countries implemented similar adjustments in response to regional inflation. These measures were followed by a gradual decline in inflation in both the US and the region, as shown in Figure \ref{fig:Eco_Signals}.

Motivated by this context, we examine the interaction of inflation signals among a group of interconnected economies, including large ones such as the United States and China, and smaller regional economies such as those in the Central American region, Panama, and the Dominican Republic (CAPARD), along with the U.S. federal funds rate, included to capture the potential influence of U.S. monetary policy on countries inflation dynamics, over the period January 2020 to November 2024. For this purpose, we will consider structured models of the form:
\[
x(t+1)=A~x(t)
\]
Here, the structure of $A$ is determined by some suitable economic interrelation constraints between the countries under consideration.

\begin{figure}[!ht]
    \centering
    \includegraphics[scale=.70]{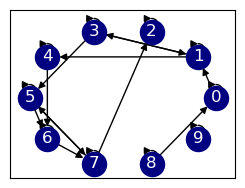}
    \caption{ Regional Inflation Network Graph}
    \label{fig:Rel_Pi_Graph}
\end{figure}

\begin{figure}[!ht]
    \centering
    \includegraphics[scale=.52]{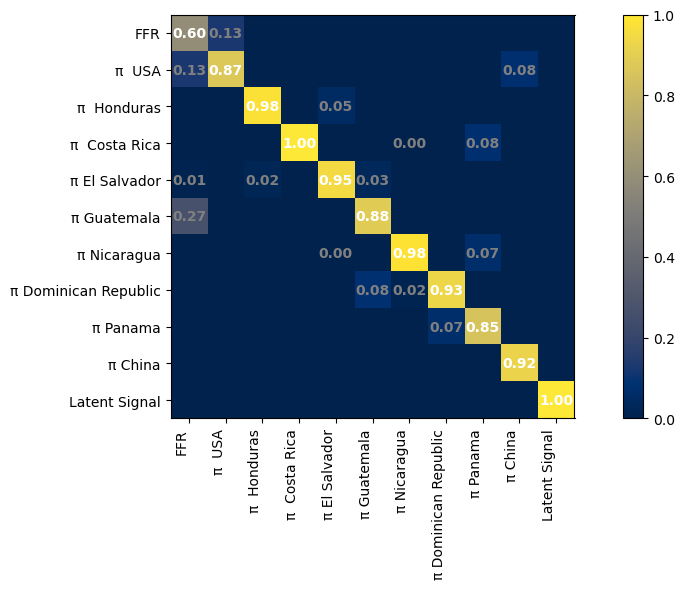}
    \caption{Contribution between the inflation networks and the US federal funds rate}
    \label{fig:contribution}
\end{figure}

Our analysis yielded notable findings. First, we identified a network of interconnections between the inflationary dynamics of the countries in the Northern Triangle of Central America inflation, and the U.S. federal funds rate, exhibiting varying lags in comparison to other countries in the region. The contributions to the inflationary states of these countries are depicted in Figure \ref{fig:contribution}, through the structural identification of the matrix and its relational behavior, as shown in the empirical relational graph presented in Figure \ref{fig:Rel_Pi_Graph}.

Although most of the inflation observed in the countries under study can be attributed to historical dynamics, we identify a direct influence of the U.S. federal funds rate on inflation in the United States, Guatemala, and El Salvador, along with indirect effects to Honduras. In addition, inflation in the United States and China has an implicit impact on these countries. Inflationary interactions among the other CAPARD countries are also evident, as illustrated in Figure \ref{fig:contribution}, suggesting the presence of an inflation transmission within the region, along with the influence of U.S. monetary policy. 
The stabilization of U.S. inflation contributed to the stabilization of inflation in Honduras during the following three periods, with faster responses in El Salvador and Guatemala.

Moreover, empirical decoupling is observed between the Federal Funds Rate, inflation in the United States and China, and inflation in the Northern Triangle countries of Central America (El Salvador, Honduras, and Guatemala). The matrix structure reveals the interactions between large and small economies in this region. In this context, the Northern Triangle countries are economies with historical trade ties to both the United States and China, as well as a significant inflow of remittances from the United States to these nations.

The behavior identified for the economic signals under consideration is illustrated in Figure \ref{fig:Eco_Signals}. The short-term predictions exhibit a high degree of accuracy compared to the inflation dynamics observed in the segment of country-level inflation rates depicted in the graph. This predictive accuracy can be attributed to the behavioral dynamics described above, particularly the country-specific Markovian structures and the mechanisms of cross-country contribution.

Second, a mathematical function has been identified that describes the dynamics of a latent signal $z(t)$, constructed from the signals studied, to stabilize the system. The model is perfectly embedded, preserving the context of inflationary dynamics and the federal funds rate, represented by the equation:

\[
x(t+1) = W_{int}x(t) + W_{lat}z(t)
\]

Here $x(t+1)$ is explained by the dynamics of the real signals in $x(t)$. The model balances interpretability and predictive capability, offering accurate predictions while preserving context.

\section{Conclusion and Future Work}

This paper presented a methodology for identifying and simulating financial and economic systems through stochastically structured reservoir computers (SSRCs). By integrating structure-preserving embeddings with graph-constrained coupling matrices, the framework yields interpretable and structurally consistent representations of dynamic processes.

Applications to a nonlinear stochastic model and a regional inflation network demonstrated the capability of SSRCs to capture nonlinear behaviors, reveal interdependencies, and reproduce emergent patterns in economic data. These results confirm the efficiency of SSRCs as modeling tools and their potential for explaining systemic change driven by agent interactions.

Future work will focus on refining latent signal modeling for stability, extending embeddings to capture shocks and long-memory effects, and developing adaptive schemes for online learning. These advances aim to consolidate SSRCs as a robust framework for analyzing complex adaptive systems.

\section*{Data Availability}

The programs and data that support the findings of this study are openly available in the DyNet-CNBS repository, reference number \cite{FVides_DyNet}.\\

\section*{Conflicts of Interest}
The authors declare that they have no conflicts of interest.

\begin{ack}
The structure preserving matrix computations needed to implement the algorithms in \S\ref{sec:results}, were performed with  {\rm Python} 3.10.4, with the support and computational resources of the National Commission of Banks and Insurance Companies ({\bf CNBS}) of Honduras. The views expressed in the article do not necessarily represent the views of the National Commission of Banks and Insurance Companies of Honduras. 
\end{ack}

\bibliography{ifacconf}             
                                                   







\end{document}